\documentclass{article}
\usepackage[left=1in,right=1in,top=1in,bottom=1in]{geometry}
\usepackage{amsmath, amssymb, amsthm, mathtools}
\usepackage{hyperref}
\usepackage{cases}
\usepackage{todonotes}
\usepackage{authblk}
\newtheorem{theorem}{Theorem}

\newtheorem{conjecture}[theorem]{Conjecture}
\newtheorem{corollary}[theorem]{Corollary}
\newtheorem{definition}[theorem]{Definition}
\newtheorem{example}[theorem]{Example}
\newtheorem{lemma}[theorem]{Lemma}
\newtheorem{note}[theorem]{Note}

\newtheorem{question}[theorem]{Question}

\begin{document}
\title{On a Conjecture about Ron Graham's Sequence}
\author[1]{Peter Kagey}
\author[2]{Krishna Rajesh}
\affil[1]{California State Polytechnic University, Pomona}
\affil[2]{Harvey Mudd College}
\date{\today}
\maketitle

\begin{abstract}
  Ron Graham's Sequence is a surprising bijection from non-negative integers to non-negative, non-prime integers that was introduced by Ron Graham in the June 1986 ``Problems'' column of \textit{Mathematics Magazine}, and which later appeared in Problem A2 of the 2013 William Lowell Putnam Mathematical Competition. We describe some properties of this function, give an algorithm for computing its values in pseudo-polynomial time, and prove a 22 year-old conjecture about an upper bound for the function.
\end{abstract}

In the June 1986 ``Problems'' column of \textit{Mathematics Magazine}, Ron Graham defined a function---which was later coined ``R. L. Graham's sequence'' in the On-Line Encyclopedia of Integer Sequences (OEIS) \cite{OEIS}---and challenged readers to prove that it was an injection \cite{Graham1986}.
The next year in the June 1987 column, a solution from Michael Reid was published. In it Reid proved not just that the function was injective, but that it was a bijection from the natural numbers to the composite numbers together with $0$ and $1$ \cite{GrahamProblemReidSolution1987}.
In 1989, the same problem appeared as Problem 4.39 in Graham, Knuth, and Patashnik's book \textit{Concrete Mathematics} \cite{GrahamKnuth1989}. In December 2013, the problem resurfaced as Problem A2 in the 2013 William Lowell Putnam Mathematical Competition \cite{PutnamBook2021}.

This paper provides a survey of what has been discovered about this sequence in the last three decades, describes the rich structure that underpins this function, and culminates in the resolution of a 2002 conjecture about an upper bound for the function.

In Section \ref{sec:Preliminaries}, we define Ron Graham's sequence, give bounds on the function, and describe Michael Reid's bijectivity proof.
In Section \ref{sec:LinearAlgebra} we discuss methods for computing the sequence by encoding it as a system of linear equations over the finite field of two elements, $\mathbb{F}_2$.
In Section \ref{sec:RobertWilsonConjecture}, we give a proof of a 2002 conjecture from Robert G. Wilson v which appeared in co-author Krishna Rajesh's undergraduate thesis.
Finally, in Section \ref{sec:Conjectures}, we give some natural generalizations to the problem and some related results about the algorithmic complexity of such problems.

\section{Preliminaries}
\label{sec:Preliminaries}
We start by defining our sequence, which the On-Line Encyclopedia of Integer Sequences (OEIS) catalogs as A006255 \cite{OEIS}.
\begin{definition}
  Let \textbf{Ron Graham's sequence} be the sequence $\left(g(n)\right)_{n=0}^\infty$ where \(g(n) = k\) is the least integer $k$ such that there exists an increasing sequence \[
    n = a_1 < a_2 < \cdots < a_t = k
  \] such that the product \(
    a_1 a_2 \cdots a_t
  \) is square.
\end{definition}
It is helpful to illustrate this definition with examples.
\begin{example}
The following table gives values of $g(n)$ for $0 \leq n \leq 15$, along with corresponding increasing sequences that begin with $n$, end with $g(n)$, and whose product is square.
\begin{alignat*}{5}
  & g(0)  &&= 0  &\qquad \text{with} \qquad && 0                                      & = 0^2    \\
  & g(1)  &&= 1  &\qquad \text{with} \qquad && 1                                      & = 1^2    \\
  & g(2)  &&= 6  &\qquad \text{with} \qquad && 2 \cdot 3 \cdot 6                      & = 6^2    \\
  & g(3)  &&= 8  &\qquad \text{with} \qquad && 3 \cdot 6 \cdot 8                      & = 12^2   \\
  & g(4)  &&= 4  &\qquad \text{with} \qquad && 4                                      & = 2^2    \\
  & g(5)  &&= 10 &\qquad \text{with} \qquad && 5 \cdot 8 \cdot 10                     & = 20^2   \\
  & g(6)  &&= 12 &\qquad \text{with} \qquad && 6 \cdot 8 \cdot 12                     & = 24^2   \\
  & g(7)  &&= 14 &\qquad \text{with} \qquad && 7 \cdot 8 \cdot 14                     & = 28^2   \\
  & g(8)  &&= 15 &\qquad \text{with} \qquad && 8 \cdot 10 \cdot 12 \cdot 15           & = 120^2  \\
  & g(9)  &&= 9  &\qquad \text{with} \qquad && 9                                      & = 3^2    \\
  & g(10) &&= 18 &\qquad \text{with} \qquad && 10 \cdot 12 \cdot 15 \cdot 18          & = 180^2  \\
  & g(11) &&= 22 &\qquad \text{with} \qquad && 11 \cdot 18 \cdot 22                   & = 66^2   \\
  & g(12) &&= 20 &\qquad \text{with} \qquad && 12 \cdot 15 \cdot 20                   & = 60^2   \\
  & g(13) &&= 26 &\qquad \text{with} \qquad && 13 \cdot 18 \cdot 26                   & = 78^2   \\
  & g(14) &&= 21 &\qquad \text{with} \qquad && 14 \cdot 15 \cdot 18 \cdot 20 \cdot 21 & = 1260^2 \\
  & g(15) &&= 24 &\qquad \text{with} \qquad && 15 \cdot 18 \cdot 20 \cdot 24          & = 360^2
\end{alignat*}
\label{ex:initialTerms}
\end{example}
The existence of an increasing sequence $a_1 < a_2 < \cdots < a_t$ with a square product gives an upper bound for $g(n)$, namely $g(n) \leq a_t$. In particular, $g(n)$ is the least upper bound over all such sequences.
\begin{example} 
  We can see that $g(2) \leq 6$ because $2 \cdot 3 \cdot 6$ is a perfect square. To check that $g(2) = 6$, we must also prove that there is no increasing sequence that starts with $2$, ends with $2$, $3$, $4$, or $5$, and has a product that is square. We can check this explicitly, since the following eight products are all non-square: \[
    \begin{matrix*}[l]
      2 = 2  \qquad & 
      2 \cdot 3 = 6  \qquad & 
      2 \cdot 4 = 8    \qquad & 
      2 \cdot 3 \cdot 4 = 24
      \\
      2\cdot5 = 10  \qquad &
      2\cdot3\cdot5 = 30  \qquad &
      2\cdot4\cdot5 = 40  \qquad &
      2\cdot3\cdot4\cdot5 = 120.
    \end{matrix*}
  \]
  \label{ex:bruteforceCheck}
\end{example}
We describe a technique in Section \ref{sec:LinearAlgebra} that allows us to compute $g(n)$ without explicitly checking each of the $2^{g(n)-n-1}$ smaller sequences.

\subsection{Square-product and corresponding sequences}
We will be referring to increasing sequences whose product is square frequently, so it is useful to give them a name.
\begin{definition}
    A \textbf{square-product sequence} or a sequence with the \textbf{square-product property} is an increasing sequence $a_1 < a_2 < \cdots < a_t$ such that $a_1 a_2 \cdots a_t$ is a square number.
\end{definition}
In particular, Ron Graham's sequence gives the least $a_t$ over all square-product sequences that begin with $n$.

\begin{definition}
  A \textbf{corresponding sequence for} $\boldsymbol{g(n)}$ is a square-product sequence that starts with $n$ and ends with $g(n)$.
\end{definition}

Note that the definition of $g(n)$ only considers the largest term of the sequence, so $2 < 3 < 6$ and $2 < 3 < 4 < 6$ are both corresponding sequences for $g(2)$. In Subsection \ref{subsec:countingCorrespondingSequences}, we will prove a perhaps surprising fact: for each $n$ there are exactly $2^N$ distinct corresponding sequences for $g(n)$ for some positive integer $N$.

\subsection{An upper bound}
The definition of Ron Graham's sequence presumes the existence of an increasing sequence starting with $n$ with a square product. We know that we can always find such a sequence because $n\cdot4n = (2n)^2$ for any integer $n$. This gives us the upper bound $g(n) \leq 4n$, but we can prove an even tighter bound on $g(n)$.

\begin{lemma}
  Let $n \geq 4$, then we can bound $g$ above by $g(n) \leq 2n$.
  \label{lem:gn<2n}
\end{lemma}
\begin{proof}
  We will find an integer $k$ such that $n < 2k^2 < 2n$, and therefore \(n \cdot 2k^2 \cdot 2n = (2kn)^2\).

  Note that $n/2 < k^2 < n$, and thus $k$ is in the interval $(\sqrt{n/2}, \sqrt{n})$. If an interval has length greater than $1$, it must contain an integer. The size of this interval is \[
    \sqrt{n} - \sqrt\frac{n}{2} = \left(1 - \frac{1}{\sqrt{2}}\right)\sqrt{n}.
  \] Setting this to be greater than $1$ and solving for $n$ yields \[
    n > \frac{1}{\left(1-\frac{1}{\sqrt{2}}\right)^2} = 6 + 4 \sqrt2 \approx 11.7
  \]
  Thus there exists such a $k$ whenever $n \geq 12$.
  Then it is enough to check that $g(n) \leq 2n$ for $4 \leq n \leq 11$, which is done in Example \ref{ex:initialTerms}.
\end{proof}

We will see in Subsection \ref{subsec:Algorithm} that good upper bounds are useful in practically computing Ron Graham's sequence. Although we will discuss a stronger upper bound in Note \ref{note:betterUpperBound}. In particular, this bound is tight for all primes greater than $3$.

\begin{lemma}
    Let $p > 3$ be prime, then $g(p) = 2p$.
    \label{lem:gn<2nIsTight}
\end{lemma}
\begin{proof}
    Let $a_1 < a_2 < \cdots < a_t$ be a corresponding sequence to $g(p)$ so that $a_1 = p$. Since $p \mid a_1a_2 \cdots a_t$ and the sequence has the square-product property, $p^2$ must divide the product. Thus the sequence must contain another term that has $p$ as a factor. Since $2p$ is the smallest integer greater than $p$ that has $p$ as a factor, the corresponding sequence must contain a term that is greater than or equal to $2p$. Thus $a_t \geq 2p$, and so $g(p) \geq 2p$. 
    
    Since by Lemma \ref{lem:gn<2n} $g(n) \leq 2n$, we can see that $g(p) = 2p$ for all primes $p > 3$.
\end{proof}

\subsection{Injectivity of Ron Graham's Sequence}

Now, we will solve the problem from Ron Graham that appeared in the June 1986 ``Problems'' column of \textit{Mathematics Magazine} and again as Problem A2 on the 2013 Putnam exam: proving that $g$ is injective. We will also prove that $g$ is surjective onto the non-prime non-negative integers. In both cases, our proof will loosely follow proof Michael Reid's solution to this problem appeared in the June 1987 column \cite{GrahamProblemReidSolution1987}.

\begin{theorem}
  Ron Graham's sequence $g \colon \mathbb{N}_{\geq 0} \to \mathbb{N}_{\geq 0}$ is injective.
  \label{thm:GrahamInjective}
\end{theorem}
\begin{proof}
  Suppose that $g$ is not injective, so that there is some $g(n) = g(n')$
  for $n < n'$.
  Then there exist corresponding sequences for $g(n)$ and $g(n')$ given by \[
    n = a_1 < a_2 < \cdots < a_t = g(n)
  \] and  \[
    n' = a'_1 < a'_2 < \cdots < a'_{t'} = g(n')
  \] such that $a_1a_2 \cdots a_t$ and $a'_1a'_2 \cdots a'_{t'}$ are both square, and
  thus $a_1a_2 \cdots a_ta'_1a'_2 \cdots a'_{t'}$ is square.

  By removing all terms that appear in both sequences (including $a_t$ and $a'_{t'}$),
  we get a resulting sequence which is also square.
  Thus $g(n) \leq \max(a_{t-1}, a'_{t'-1}) < a_t = a'_t = g(n')$, a contradiction.
\end{proof}

\subsection{Surjectivity of Ron Graham's Sequence}
In this section, we will show that $g$ is surjective onto the non-negative, non-prime integers---that is, the composite numbers together with $0$ and $1$. We begin by proving that $g(n)$ is never prime.
\begin{lemma}
  Let $\mathbb{P}$ denote the set of prime numbers. Then $g(n) \in (\mathbb{N}_{\geq 0} \setminus \mathbb{P})$ for all non-negative integers $n$.
\end{lemma}
\begin{proof}
  When $n = 0$, $g(0) = 0 \in (\mathbb{N}_{\geq 0} \setminus \mathbb{P})$. For $n > 0$ and a square-product sequence \(
    n = a_1 < a_2 < \cdots < a_t
  \) we claim that $a_t$ cannot be prime. If it were, then $p^2 \mid a_1a_2\cdots a_{t-1}p$ and thus $p \mid a_1a_2\cdots a_{t-1}$. However, this cannot occur since $a_i < p$ for all $i < t$.
\end{proof}
Now, in order to prove surjectivity of $g$, we begin by defining a function, which we will later show is a right-inverse of $g$ on the set of non-prime, non-negative integers.
\begin{definition}
  Let $\bar{g}\colon \mathbb N_{\geq 0} \to \mathbb N_{\geq 0}$ be 
  $\bar{g}(k) = n$ where $n$ is the greatest integer such that there exists a square-product sequence \[
    \bar{g}(k) = n = a_1 < a_2 < \cdots < a_t = k.
  \]
\end{definition}
This is sequence in the OEIS as sequence A067565 \cite{OEIS}.
\begin{lemma}
    When $k > 0$ is non-prime, $\bar{g}(k) > 0$.
\end{lemma}
\begin{proof}
    Consider the prime factorization of $k = p_1^{e_1}p_2^{e_2}\cdot p_t^{e_t}$, and let $q_1 < q_2 < \cdots < q_t$ be the primes that appear an odd number of times in the prime factorization. Then \(q_1q_2\cdots q_tk\) is square, and so $\bar{g}(k) \geq q_1 > 0$.
\end{proof}
Now we prove that $\bar{g}$ is a right-inverse of $g$ as claimed above.
\begin{lemma}
  When $k$ is a non-prime non-negative integer, $g(\bar{g}(k)) = k$.
  \label{lem:rightInverse}
\end{lemma}
\begin{proof}
  Notice that $g(\bar{g}(k)) \leq k$ because by the definition of $\bar{g}(k)$, there exists a square-product sequence $\bar{g}(k) = a_1 < a_2 < \cdots < a_t = k$, so $k$ is an upper bound for $g(\bar{g}(k))$.

  Now we need to check that $g(\bar{g}(k))$ is not strictly less than $k$. For the sake of contradiction, suppose it were---that is, there exists a square-product sequence \(
    \bar{g}(k) = a'_1 < a'_2 < \cdots < a'_{t'} < k
  \). Since all of terms are positive, we can apply the same argument from the proof of Lemma \ref{thm:GrahamInjective} and see that removing the terms that appear in both sequences results in a sequence that starts with a number greater than $\bar{g}(k)$, ends with $k$, and whose product is square. However, this contradicts the definition of $\bar{g}(k)$, and so $g(\bar{g}(k)) = k$ whenever $k$ is non-prime.
\end{proof}
The surjectivity claim follows as a corollary.
\begin{theorem}
  Ron Graham's sequence $g \colon \mathbb N_{\geq 0} \to (\mathbb N_{\geq 0} \setminus \mathbb{P})$, is surjective onto the non-prime, non-negative integers.
  \label{thm:GrahamSurjective}
\end{theorem}
\begin{proof}
  We want to show that $g \colon \mathbb{N}_{\geq 0} \to (\mathbb{N}_{\geq 0} \setminus \mathbb{P})$ is surjective. Let $k$ be a non-negative non-prime integer, and let $n = \bar{g}(k)$. By Lemma \ref{lem:rightInverse}, $g(n) = g(\bar{g}(k)) = k$, and so $k$ is in the image of $g$.
\end{proof}
Putting this theorem together with Theorem \ref{thm:GrahamInjective} gives us the bijection we claimed at the beginning of the section.
\begin{corollary}
  Ron Graham's sequence is a bijection from the non-negative integers to the non-prime non-negative integers.
\end{corollary}

\section{The Linear Algebra of Ron Graham's Sequence}
\label{sec:LinearAlgebra}
Computing Ron Graham's sequence naively appears to take exponential time. For instance, in Example \ref{ex:bruteforceCheck}, we check that $g(2) = 6$ by providing a corresponding sequence to $g(2)$, and then computing all $2^3 = 8$ increasing sequences starting with $2$ and with maximum term strictly less than $g(2) = 6$.
Since $g(n) \leq 2n$ for $n > 3$ is a tight upper bound on $g$, this approach requires a number of checks that is exponential in $n$, and thus does not provide a practical method for computing $g(n)$ for large $n$.

For an algorithm that runs in polynomial time with respect to the size of the input (pseudo-polynomial time), we will instead make the observation that in a square-product sequence, each prime factor appears an even number of times, as shown in the following example.

\begin{example}
  The prime factorization for a corresponding sequence for $g(8)$  is
  \[
  \underbrace{2^3}_8 <
  \underbrace{2 \cdot 5}_{10} <
  \underbrace{2^2 \cdot 3}_{12} <
  \underbrace{3\cdot 5}_{15}.
\] Each prime factor appears an even number of times: $2$ appears six times and $3$ and $4$ each appear twice, thus $8 \cdot 10 \cdot 12 \cdot 15 = (2^3\cdot3\cdot5)^2 = 120^2$.
\label{ex:primeFactorMatching}
\end{example}

\subsection{An algorithm to compute \texorpdfstring{$g(n)$}{g(n)}}
\label{subsec:Algorithm}
In order to construct a sequence whose product is square, we will use a technique described by Morrison and Brillhart \cite{Morrison1975} and discussed at greater length by Pomerance in a discussion of the quadratic sieve algorithm \cite{Pomerance1996}. This algorithm exploits the fact that we only need to pair up the exponents of the prime factors of the terms---and thus by representing these exponents as vectors over the field of two elements, we find such sequences by solving a system of linear equations.

This setting up the appropriate matrix equation and solving it can be done in polynomial time with respect to the size of the input. This contrasts with the naive, brute-force approach of trying all smaller sequences, which runs in exponential time with respect to the size of the input.

\begin{definition}
  Let $n$ be an integer with prime factorization $\displaystyle n = \prod p_i^{v_i}$, where $p_i$ is the $i$-th prime.

  Then we say that the \textbf{exponent vector} $\vec{v}(n) \in \mathbb{F}_2^N$ is defined by \[
    \vec{v}(n) = \left[\begin{array}{ccccc}
      v_1 & v_2 & v_3 & \cdots & v_N
    \end{array}\right]^T, 
  \] where the entries are considered modulo $2$ and $N$ is taken to be any sufficiently large integer so that $v_i \equiv 0 \pmod 2$ for all $i > N$. 
\end{definition}
\begin{example}
  Expanding on Example \ref{ex:primeFactorMatching}, we can see that the exponent vectors arising from a corresponding sequence of $g(8)$ sum to the zero vector in $\mathbb{F}_2^N$: \[
    \underbrace{
      \left[\begin{array}{c}
        1 \\ 0 \\ 0 \\ 0 \\ \vdots \\ 0
      \end{array}\right]
    }_{\vec{v}(8)}
    +
    \underbrace{
      \left[\begin{array}{c}
        1 \\ 0 \\ 1 \\ 0 \\ \vdots \\ 0
      \end{array}\right]
    }_{\vec{v}(10)}
    +
    \underbrace{
      \left[\begin{array}{c}
        0 \\ 1 \\ 0 \\ 0 \\ \vdots \\ 0
      \end{array}\right]
    }_{\vec{v}(12)}
    +
    \underbrace{
      \left[\begin{array}{c}
        0 \\ 1 \\ 1 \\ 0 \\ \vdots \\ 0
      \end{array}\right]
    }_{\vec{v}(15)}
    = \underbrace{
      \left[\begin{array}{c}
        0 \\ 0 \\ 0 \\ 0 \\ \vdots \\ 0
      \end{array}\right]
    }_{\vec{0}}.
  \]
\end{example}
\begin{lemma}
  \label{lem:squareProductWhenVectorsSum}
  The sequence $a_1 < a_2 < \cdots < a_n$ is a square-product sequence if and only if 
  $\vec{v}(a_1) + \vec{v}(a_2) + \cdots + \vec{v}(a_n) = \vec{0}$.
\end{lemma}
\begin{proof} 
  Notice that for all integers $x, y \in \mathbb{N}_\geq{0}$, $\vec{v}(xy) = \vec{v}(x) + \vec{v}(y)$. Moreover, an integer has an exponent vector equal to the zero vector only if and only if all of its prime factors appear an even number of times.
\end{proof}

\begin{definition}
  Let $A(\ell, r)$ be the matrix whose columns are given by $\vec{v}(\ell), \vec{v}(\ell + 1), \dots, \vec{v}(r)$, where the number of rows in the matrix is chosen to be at least $\pi(r)$, where $\pi$ is the prime counting function.
\end{definition}
\begin{lemma}
  \label{lem:matrixEquation}
  When $n$ is non-square there exists a square-product sequence starting with $n$ that is bounded above by $r$ whenever there exists a particular solution to the matrix equation \(
    A(n+1, r)\vec{b} + \vec{v}(n) = \vec 0
  \) over the field $\mathbb{F}_2$, or equivalently, \[
    A(n+1, r)\vec{b} = \vec{v}(n).
  \]
\end{lemma}
\begin{proof}
  A particular solution to this matrix equation corresponds to a linear combination of the vectors that together with $\vec{v}(n)$ sum to $\vec{0}$: \[
    c_{n+1}\vec{v}(n+1) +
    c_{n+2}\vec{v}(n+2) +
    \cdots
    + c_{r}\vec{v}(r) 
    + \vec{v}(n) = \vec{0}.
  \] Since this is over the finite field $\vec{F}_2$, each $c_i \in \{0, 1\}$, and so the linear combination of the exponent vectors corresponds to a subset of $\{n+1, n+2, \dots, r\}$ By Lemma \ref{lem:squareProductWhenVectorsSum}, the product of this subset is square when multiplied by $n$, and this gives a square-product sequence that begins with $n$ and ends with a number that is less than or equal to $r$.
\end{proof}
\begin{example}
  Let $n = 8$, and consider $A(9,15)\vec{b} = \vec{v}(8)$: \[
    \left[\begin{array}{cccc}
      \vline & \vline &  & \vline \\
      \vec{v}(9) & \vec{v}(10) & \cdots & \vec{v}(15) \\
      \vline & \vline & & \vline
    \end{array}\right]
    \vec{b} = \vec{v}(8)
  \]
  Specifically \[
    \left[\begin{array}{ccccccc}
      0 & 1 & 0 & 0 & 0 & 1 & 0 \\ 
      0 & 0 & 0 & 1 & 0 & 0 & 1 \\ 
      0 & 1 & 0 & 0 & 0 & 0 & 1 \\ 
      0 & 0 & 0 & 0 & 0 & 1 & 0 \\ 
      0 & 0 & 1 & 0 & 0 & 0 & 0 \\ 
      0 & 0 & 0 & 0 & 1 & 0 & 0    
    \end{array}\right]\vec{b} = \left[\begin{array}{c}
      1 \\ 
      0 \\ 
      0 \\ 
      0 \\ 
      0 \\ 
      0    
    \end{array}
    \right],
  \] which has a particular solution of \[
    \vec{b} = [0,1,0,1,0,0,1]^T,
  \] where the $1$s correspond to columns $\vec{v}(10)$, $\vec{v}(12)$, and $\vec{v}(15)$.

  Thus $8 < 10 < 12 < 15$ is an increasing sequence whose product is square.
  Since $A(9,14)\vec{b} = \vec{v}(8)$ does not have any solutions, this means that $g(8) = 15$, and the above sequence is a corresponding sequence to $g(8)$.
\end{example}
\subsection{Counting corresponding sequences}
\label{subsec:countingCorrespondingSequences}
This perspective of computing $g$ by solving matrix equations over $\mathbb{F}_2$ has other implications as well, including being able count the number of corresponding sequences for $g(n)$.
\begin{lemma}
  When $n$ is non-square, the set of corresponding sequences for $g(n)$ is in bijection with the null space of $A(n+1, g(n))$.
\end{lemma}
\begin{proof}
    Each particular solution to a matrix equation of this form corresponds to a square-product sequence, and the set of particular solutions is in bijection with the null space. Moreover, by construction, each particular solution corresponds to a sequence whose maximal value is $g(n)$.
\end{proof}
Because we are solving matrix equations over a finite field, there are a finite number of vectors in the null space.
\begin{corollary}
  Let $N$ be the dimension of the null space of $A(n+1, g(n))$: $\operatorname{nullity}\!\left(A\left(n+1,g(n)\right)\right) = N$. Then there are $2^N$ corresponding sequences for $g(n)$.
\end{corollary}
\begin{proof}
  Any linear combination of the basis of the null space gives an element of the null space, and there are $2$ choices for each coefficient, resulting in a total of $2^N$ elements of the null space and thus $2^N$ corresponding sequences for $g(n)$.
\end{proof}
Thus given the null space of such a matrix, we can systematically enumerate all of the corresponding sequences for $g(n)$.
\begin{example}
  Let $n = 11$ and note that $11 < 18 < 22$ is a corresponding sequence for $g(11) = 22$.
  The following four vectors form a basis of the nullspace of $A(11+1,g(11)) = A(12,22)$:
  \[
    \begin{array}{ccccc cccccc}
      12 & 13 & 14 & 15 & 16 & 17 & 18 & 19 & 20 & 21 & 22 \\
      \hline
      1 &  0 &  1 &  0 &  0 &  0 &  1 &  0 &  0 &  1 &  0 \\
      1 &  0 &  0 &  1 &  0 &  0 &  0 &  0 &  1 &  0 &  0 \\
      0 &  0 &  0 &  0 &  1 &  0 &  0 &  0 &  0 &  0 &  0
    \end{array}
  \]
  These basis vectors correspond to the products 
  $12\cdot14\cdot18\cdot21 = 252^2$, 
  $12\cdot15\cdot20 = 60^2$, and
  $16 = 4^2$.
  So there are $2^3 = 8$ corresponding sequences for $g(13)$:
  \begin{alignat*}{2}
    & 11 \cdot 18 \cdot 22                                     &&= 66^2    \\
    & 11 \cdot 12 \cdot 14 \cdot 21 \cdot 22                   &&= 924^2   \\
    & 11 \cdot 12 \cdot 15 \cdot 18 \cdot 20 \cdot 22          &&= 3960^2  \\
    & 11 \cdot 14 \cdot 15 \cdot 20 \cdot 21 \cdot 22          &&= 4620^2  \\
    & 11 \cdot 16 \cdot 18 \cdot 22                            &&= 264^2   \\
    & 11 \cdot 12 \cdot 14 \cdot 16 \cdot 21 \cdot 22          &&= 3696^2  \\
    & 11 \cdot 12 \cdot 15 \cdot 16 \cdot 18 \cdot 20 \cdot 22 &&= 15840^2 \\
    & 11 \cdot 14 \cdot 15 \cdot 16 \cdot 20 \cdot 21 \cdot 22 &&= 18480^2
  \end{alignat*}
\end{example}

\begin{example}
OEIS sequence A260510 is defined as $\operatorname{A260510}(n) = N$ where there are $2^N$ corresponding sequences for $g(n)$ \cite{OEIS}. 
This sequence begins as follows. \[
    \begin{array}{r|cccccccccccccccccccc}
        n & 1 &  2 &  3 &  4 &  5 &  6 &  7 &  8 &  9 &
        10 & 11 & 12 & 13 & 14 & 15 & 16 & 17 & 18  \\
        \hline
        \operatorname{A260510}(n) & 0 &  1 &  1 &  0 &  1 &  1 &  1 &  1 &  0 &
        1 &  3 &  1 &  4 &  1 &  1 &  0 &  6 &  1
    \end{array}
  \]
  OEIS sequence A259527 is defined as $\operatorname{A259527}(n) = 2^N$ where there are $2^N$ corresponding sequences for $g(n)$.
  This sequence begins as follows. \[
    \begin{array}{r|cccccccccccccccccccc}
        n & 1 &  2 &  3 &  4 &  5 &  6 &  7 &  8 &  9 &
        10 & 11 & 12 & 13 & 14 & 15 & 16 & 17 & 18  \\
        \hline
        \operatorname{A259527}(n) & 1 & 2 & 2 & 1 & 2 & 2 & 2 & 2 & 1 & 2 & 8 & 2 & 16 & 2 & 2 & 1 & 64 & 2
    \end{array}
  \]
\end{example}

\section{A conjecture of Robert G. Wilson v} 
\label{sec:RobertWilsonConjecture}

In January 2002, Robert G. Wilson v made an observation in the OEIS that corresponding sequences for $g(n)$ never appear to have length $2$. In this subsection, we prove this claim by explicitly constructing an appropriate sequence.

We begin by defining the minimal length of a corresponding sequence, which appears in the OEIS as A066400 \cite{OEIS}.
\begin{definition}
  Let $T(n)$ be the least $t$ such that there exists a corresponding sequence for $g(n)$ of length $t$, \[
    n = a_1 < a_2 < \cdots < a_t = g(n).
  \]
  \label{def:Tn}
\end{definition}

For all $t \leq 14$ except for $t = 2$, we have examples of integers $n$ such that $T(n) = t$.
\begin{example}
  The least $n$ such that $T(n) = t$ is given in the following table \[
    \begin{array}{l|llllllllllllll}
    T(n) & 1 & 2 & 3 & 4 &  5 &  6 &  7 &   8 &   9 &  10 &   11 &   12 &   13 &   14 \\ \hline
    n & 1 & - & 2 & 8 & 14 & 52 & 99 & 589 & 594 & 595 & 1566 & 1961 & 3465 & 5301
    \end{array}
  \]
  In particular, if \[
    5301 = a_1 < a_2 < \cdots < a_t = g(5301) = 5375,
  \]
  is a square-product sequence (and thus a corresponding sequence for $g(5301)$) then $t \geq 14$.
  \label{ex:RecordT}
\end{example}

In order to prove the conjecture, we define a related function, which is in the OEIS as A072905 \cite{OEIS}.
\begin{definition}
  Let $f(n)$ be the least integer $k > n$ such that $nk$ is square.
\end{definition}
\begin{lemma}
  $T(n) = 2$ if and only if $g(n) = f(n)$.
\end{lemma}
\begin{proof} First, suppose that $T(n) = 2$, then $n = a_1 < a_2 = g(n)$ and $a_1a_2 = ng(n)$ is square---and moreover, $g(n)$ is the least integer with this property. Thus $g(n) = f(n)$ by the definition of $f(n)$.

~

Next, suppose that $g(n) = f(n)$. Then $T(n) \leq 2$ because $nf(n)$ is square, and $T(n) > 1$ because $f(n) > n$ by definition.
\end{proof}

In order to show that $g(n) \neq f(n)$, it is useful to factor the largest square from $n$.
\begin{definition}
  Let $n = mr^2$ where $m$ is the largest integer such that $n/m$ is square.
  Then $m$ is called the \textbf{squarefree part of $n$} and $r^2$ is called the \textbf{largest square dividing $n$}.
\end{definition}

\begin{lemma}
  Let $n = mr^2$ where $m$ is the squarefree part of $n$. 
  Then $f(n) = m(r+1)^2$.
  \label{lem:t2iffF=G}
\end{lemma}
\begin{proof}
  Notice that in order for $nk$ to be square, $n$ and $k$ must have the same squarefree part, $m$, and so $f(n) = ma^2$ for some integer $a$. Since $f(n) > n$, $a > r$, and so $a = r+1$ is the smallest choice for $a$.
\end{proof}
The following theorem was first proven in the senior thesis of the co-author Krishna Rajesh. \cite{Rajesh2024}
\begin{theorem}[Robert G. Wilson's conjecture]
  \label{thm:RGWConjecture}
  For all non-negative integers $n$, the length of the shortest corresponding sequence to $g(n)$, $T(n)$, is not equal to $2$.
\end{theorem}
\begin{proof}
  First, notice that if $n$ is square, then $T(n) = 1 \neq 2$. Therefore assume $n$ is not a square, and let $n = mr^2$ where $m > 1$ is the squarefree part of $n$.

  Next, let $s = mr + 1$, and consider the increasing sequence \[
    \underbrace{mr^2}_n < \underbrace{mr^2 + r}_{rs} < \underbrace{mr^2 + mr}_{mr(r+1)} < \underbrace{mr^2 + mr + r + 1}_{(r+1)s}
  \] whose product is $(mr^2(r+1)s)^2$. Thus $g(n) \leq (r+1)s$.

  Since $m > 1$, we can directly check that \[
    g(n) \leq (r+1)s = mr^2 + mr + r + 1 < m(r+1)^2 = mr^2 + 2mr + 1 = f(n),
  \] 
  By Lemma \ref{lem:t2iffF=G}, since $g(n) \neq f(n)$, we conclude that $T(n) \neq 2$.
\end{proof}

\begin{corollary}
  No corresponding sequences have length $2$.
\end{corollary}

\begin{note}
  \label{note:betterUpperBound}
  When $n = mr^2 > 4$ is not squarefree, $g(n) \leq mr^2 + mr + r + 1$ is a tighter upper bound than the bound $g(n) < 2n$ from Lemma \ref{lem:gn<2n},
  which allows for a smaller matrix computation in Lemma \ref{lem:matrixEquation}.
\end{note}

\section{Conjectures and related questions}
\label{sec:Conjectures}
We finish with some conjectures and related questions.
\subsection{Conjectures}
Recall that Example \ref{ex:RecordT} gives examples where $T(n) = t$ for all integers $t \leq 14$, besides $t = 2$. We have checked by computer that $T(n) < 15$ for all $n \leq 10000$.

\begin{conjecture}
  The function $T(n)$ which gives the minimal length of a corresponding sequence is unbounded. Moreover, it is surjective onto the set of positive integers not equal to $2$.
\end{conjecture}
We make an even stronger conjecture!
\begin{conjecture}
  For each positive integer $t \neq 2$, the set \(
    \{n \in \mathbb{N}_{\geq 0} \mid T(n) = t\}
  \) is infinite.
\end{conjecture}

In Lemma \ref{lem:gn<2n}, we saw that $g(n) \leq 2n$, and in Lemma \ref{lem:gn<2nIsTight}, we saw that this bound is tight, we conjecture that the bound is (almost) only achieved on prime inputs.
\begin{conjecture}
   We achieve the upper bound $g(n) = 2n$ if and only if $n$ is prime or $n = 6$.
\end{conjecture}

\subsection{Related questions and generalizations}
\begin{definition}
    Let a corresponding sequence be called \textbf{primitive} if has no proper non-empty subsequences with the square-product property. 
\end{definition}
The number of primitive corresponding sequences is counted by OEIS sequence A291634 \cite{OEIS}.
\begin{question}
    Is there a psuedo-polynomial time algorithm to compute the number of primitive corresponding sequences?
\end{question}
The techniques in Section \ref{sec:LinearAlgebra} suggest generalizations over other finite fields, which suggests the following generalized definition.
\begin{definition}
  Say that $g_m(n)$ is the least integer $k$ such that there exists a sequence 
  \[
    n = a_1 \leq a_2 \leq \cdots \leq a_t = k,
  \] such that $a_1a_2 \cdots a_t = N^m$ for some integer $N$, and no integer appears in the sequence $m$ or more times.
  \label{def:generalizedSequence}
\end{definition}
Now we can ask all of the analogous questions about $g_m(n)$:
\begin{question}
  Is there a polynomial time algorithm for computing $g_m(n)$ when $m$ is not prime?
\end{question}
\begin{question}
  We saw that $g(n) = g_2(n)$ is surjective onto the non-prime non-negative integers. Does the image of $g_m(n)$ have any notable properties?
\end{question}
\begin{question}
  There were $2^N$ corresponding sequences for $g_2(n)$. Does the set of corresponding sequences to $g_m(n)$ have any notable properties when $m$ is not prime?
\end{question}
\begin{definition}
  Let $T_m(n)$ and $T'_m(n)$ be the minimum number of terms and the minimum number of distinct terms over all corresponding sequences for $g_m(n)$
\end{definition}
\begin{question}
  We know that $T_2(n) = T'_2(n) \neq 2$. Do $T_m(n)$ or $T'_m(n)$ have any notable properties?
\end{question}

Instead of generalizing the power, we could also generalize the operation. 
\begin{definition}
  The OEIS defines $A300516(n) = k$ where $k$ is the least integer such that there is an increasing sequence $n = a_1 < a_2 < \cdots < a_t = k$ such that $\operatorname{LCM}(a_1, a_2, \dots, a_t)$ is square.
\end{definition}
\begin{question}
  Is there a polynomial time algorithm for computing $A300516$?
\end{question}
\begin{question}
  Does the set of corresponding sequences for $A300516(n)$ have any notable properties?
\end{question}

\subsection{Algorithms}
Lastly, it is worth mentioning a related theorem that appeared on Computer Science Stack Exchange. If in Definition \ref{def:generalizedSequence}, we require strict inequalities, then Yuval Filmus constructs a reduction to prove the following theorem.
\begin{theorem}[Yuval Filmus \cite{Filmus2020}]
  Given a finite set of integers $S$, the decision problem ``is there a non-empty subset of $S$ whose product is a cube'' is NP-hard.
\end{theorem}

\appendix

\section*{Acknowledgments}
We would like to thank Alec Jones, Art Benjamin, and Sarosh Adenwalla for their suggestions, feedback, and corrections. We acknowledge that this paper would likely not have been written without contributions from the late Reinhard Zumkeller. We would also like to extend a special thank you to Elise Lockwood---the first step of the journey only happened because you braved a trip to Oregon State University's campus in the snow in December 2013.
\bibliographystyle{plain}
\bibliography{graham.bib}

\end{document}